\documentclass[12pt]{amsart}

\usepackage[margin=3.5cm]{geometry}

\allowdisplaybreaks[4]

\usepackage{amsmath}
\usepackage{amssymb}
\usepackage{amsthm}
\usepackage{color}
\usepackage{enumerate}
\usepackage[dvips]{graphicx}
\usepackage{comment}
\usepackage{cleveref}
\usepackage{caption}

\newcommand{\T}{\mathbb{T}}

\DeclareMathOperator\SC{SC}

\def\R{\mathbb{R}}

\def\Z{\mathbb{Z}}
\def\N{\mathbb{N}}

\def\cF{\mathcal{F}}

\def\bye{\end{document}}
\def\by{\end{proof}\bye}

\def\hello{\begin{document}}
\def\fr{\frac}
\def\disp{\displaystyle}
\def\ga{\alpha}
\def\go{\omega}
\def\gep{\varepsilon}
\def\ep{\gep}
\def\mid{\,:\,}
\def\gb{\beta}
\def\gam{\gamma}
\def\gd{\delta}
\def\gz{\zeta}
\def\gth{\theta}
\def\gk{\kappa}
\def\gl{\lambda}
\def\gL{\Lambda}
\def\gs{\sigma}
\def\gf{\varphi}
\def\tim{\times}
\def\aln{&\,}
\def\ol{\overline}
\def\ul{\underline}
\def\pl{\partial}
\def\hb{\text}
\def\cF{\mathcal{F}}
\def\Int{\mathop{\text{int}}}

\def\gG{\varGamma}
\def\lan{\langle}
\def\ran{\rangle}
\def\cD{\mathcal{D}}
\def\cB{\mathcal{B}}
\def\bcases{\begin{cases}}
\def\ecases{\end{cases}}
\def\balns{\begin{align*}}
\def\ealns{\end{align*}}
\def\balnd{\begin{aligned}}
\def\ealnd{\end{aligned}}
\def\bgat{\begin{gathered}}
\def\egat{\end{gathered}}
\def\1{\mathbf{1}}
\def\bproof{\begin{proof}}

\def\eproof{\end{proof}}

\theoremstyle{definition}
\newtheorem{definition}{Definition}
\theoremstyle{plain}
\newtheorem{theorem}[definition]{Theorem}
\newtheorem{corollary}[definition]{Corollary}
\newtheorem{lemma}[definition]{Lemma}
\newtheorem{proposition}[definition]{Proposition}
\theoremstyle{remark}
\newtheorem{remark}[definition]{Remark}
\newtheorem{notation}[definition]{Notation}

\def\red#1{\textcolor{red}{#1}}
\def\blu#1{\textcolor{blue}{#1}}
\def\beq{\begin{equation}}
\def\eeq{\end{equation}}
\def\bthm{\begin{theorem}}
\def\ethm{\end{theorem}}
\def\bproof{\begin{proof}}
\def\eproof{\end{proof}}

\usepackage[abbrev]{amsrefs}
\def\eqr#1{\eqref{#1}}
\def\bmat{\begin{pmatrix}}
\def\emat{\end{pmatrix}}

\newcommand{\Pmo}{\mathcal{P}^-_{1}}
\newcommand{\Ppo}{\mathcal{P}^+_{1}}
\newcommand{\Pmk}{\mathcal{P}^-_{k}}
\newcommand{\Ppk}{\mathcal{P}^+_{k}}
\def\diag{\operatorname{diag}}
\def\rT{\mathrm{T}}
\newcommand{\Rn}{{\mathbb R}^N}
\def\IN{\text{ in } }

\def\AND{\text{ and }}
\def\FOR{\text{ for }}
\def\FORALL{\text{ for all }}
\def\ON{\text{ on }} \def\OR{\text{ or }} 
\def\OW{\text{ otherwise}}
\def\IF{\hb{ if }}
\def\WITH{\text{ with }}

\def\rmb{\mathrm{b}}
\def\I{\mathbb{I}}
\def\du#1{\left\lan#1\right\ran}
\def\bald{\begin{aligned}}
\def\eald{\end{aligned}}
\def\stm{\setminus}
\def\t{\tau}
\def\SC-{\operatorname{SC}^-}
\def\lip{\operatorname{lip}}
\def\B{\operatorname{\mathbb{B}}}

\def\Sp{\operatorname{Sp}}
\def\P{\operatorname{\mathbb{P}}}
\def\erf{\eqref}
\def\cG{\mathcal{G}}
\def\cE{\mathcal{E}}
\def\pr{\,^\prime}
\def\gX{\Xi}
\def\gS{\Sigma}
\def\cV{\mathcal{V}}
\def\cW{\mathcal{W}}
\def\Hall{(H_i)_{i\in\I}}\def\Lall{(L_i)_{i\in\I}}

\def\fC{\frak{C}}
\def\fM{\frak{M}}
\def\0{\mathbf{0}}

\renewcommand{\subjclassname}{%
\textup{2010} Mathematics Subject Classification}

\title[Vanishing discount problem]{The vanishing discount problem for 
monotone systems of Hamilton-Jacobi equations: \\ a counterexample to the full convergence$^\sharp$}

\thanks{$^\sharp$This paper is dedicated to Neil Trudinger on the occasion of his 80th birthday. 
\author[H. Ishii]{Hitoshi Ishii}
\address[\textsc{Hitoshi Ishii}]{Institute for Mathematics and Computer Science\newline
\indent Tsuda University  \newline
 \indent   2-1-1 Tsuda, Kodaira, Tokyo, 187-8577 Japan.
}
\email{hitoshi.ishii@waseda.jp} }


\def\alert#1{\begin{color}{red}#1 \end{color}}

\begin{document}
\maketitle
\begin{abstract} 
In recent years there has been intense interest in the vanishing discount problem for Hamilton-Jacobi equations. 
In the case of the scalar equation, B. Ziliotto has recently given an example of the Hamilton-Jacobi equation having non-convex Hamiltonian in the gradient variable, for which the full convergence of the solutions does not hold as the discount factor tends to zero.  
We give here an explicit example of nonlinear monotone systems of Hamilton-Jacobi equations having convex  Hamiltonians in the gradient variable, for which the full convergence of the solutions fails as the discount factor goes to zero. 
\end{abstract}

\tableofcontents

\def\rmc{\mathrm{c}}

\section{Introduction}
We consider the system of Hamilton-Jacobi equations
\beq \label{eq:1}
\bcases
\gl u_1(x)+H_1(Du_1(x))+B_1(u_1(x),u_2(x))=0 \ & \IN \T^n, \\
\gl u_2(x)+H_2(Du_2(x))+B_2(u_1(x),u_2(x))=0 \ & \IN \T^n,
\ecases
\eeq
where $\gl>0$ is a given constant, the functions $H_i : \R^n\to \R$ and $B_i : \R^2 \to R$, with $i=1,2$, are given continuous functions, and $\T^n$ denotes the $n$-dimensional flat torus $\R^n/\Z^n$.

In a recent paper \cite{IJ}, the authors have investigated the vanishing discount problem for a nonlinear monotone system of Hamilton-Jacobi equations 
\beq\label{eq:1.2}\bcases
\gl u_1(x)+G_1(x,Du_1(x),u_1(x),u_2(x),\ldots,u_m(x))=0 \ & \IN \T^n, \\
\phantom{\gl u_1(x)+G_1(x,Du_1(x),u_1(x),u_2(x)}\vdots &\\
\gl u_m(x)+G_m(x,Du_m(x), u_1(x),u_2(x),\ldots,u_m(x))=0 \ & \IN \T^n,
\ecases
\eeq
and established under some hypotheses on the $G_i 
\in C(\T^n\tim\R^n\tim \R^m)$ that, when $u^\gl=(u_1^\gl,\ldots,u_m^\gl)\in C(\T^n)^m$ denoting the (viscosity) solution 
of \erf{eq:1.2}, the whole family $\{u^\gl\}_{\gl>0}$ converges in $C(\T^n)^m$ 
to some $u^0\in C(\T^n)^m$ as $\gl\to 0+$.  The constant $\gl>0$ in the above system is the so-called \emph{discount factor.} 

The hypotheses on the system are the convexity, coercivity, and monotonicity of the $G_i$ as well as the solvability of \erf{eq:1.2}, with $\gl=0$.  Here the convexity 
of $G_i$ is meant that the functions $\R^n\tim\R^m\ni (p,u)\mapsto 
G_i(x,p,u)$ are convex. We refer to \cite{IJ} for the precise 
statement of the hypotheses.  

Prior to work \cite{IJ}, there have been many contributions to the question about 
the whole family convergence (in other words, the full convergence) under the vanishing discount, which we refer to 
\cites{IJ, DZ2, DFIZ, IMT1, IMT2, IS, CCIZ} and the references therein. 

In the case of the scalar equation, B. Ziliotto \cite{Zi} has recently shown an example of the Hamilton-Jacobi equation having non-convex Hamiltonian in the gradient variable for which the full convergence does not hold.  
In Ziliotto's approach, the first step is to find a system of two algebraic 
equations
\beq\label{two-sys}\left\{\bald
&\gl u+f(u-v)=0,
\\&
\gl v+g(v-u)=0,
\eald\right.
\eeq
with two unknowns $u,v\in\R$ and with a parameter $\gl>0$ as the discount factor, for which the solutions $(u^\gl,v^\gl)$ stay bounded and fail to fully converge as $\gl \to 0+$. Here, an ``algebraic'' equation means not to be a functional equation. 
The second step is to 
interporate the two values $u^\gl$ and $v^\gl$ to get a function of $x\in\T^1$ 
which satisfies a scalar non-convex Hamilton-Jacobi equation in $\T^1$.  

In the first step above, Ziliotto constructs $f,g$ based on a game-theoretical 
and computational argument, and the formula for $f, g$ is of 
the minimax type and not quite explicit.   In \cite{IH}, the author has 
reexamined the system given by Ziliotto, with a slight generality, as a 
counterexample for the full convergence in the vanishing discount.  

Our purpose in this paper is to present a system \erf{two-sys}, with an explicit 
formula for $f,g$,  for which the solution $(u^\gl,v^\gl)$ does not fully 
converge to a single point in $\R^2$.  A straightforward consequence is that 
\erf{eq:1}, with $
B_1(u_1,u_2)=f(u_1-u_2)$ and $B_2(u_1,u_2)=g(u_2-u_1)$,
has a solution given by 
\[
(u_1^\gl(x),u_2^\gl(x))=(u,v) \ \ \FOR x\in \T^n,
\]
under the assumption that $H_i(x,0)=0$ for all $x\in\T^n$, and therefore, gives an example of a discoutned system of Hamilton-Jacobi equations, the solution of which fails to satisfy the full convergence as the discount factor goes to zero.

The paper consists of two sections. This introduction is followed 
by Section 2, the final section, which is divided into three subsections.   
The main results are stated in the first subsection of Section 2, 
the functions $f,g$, the key elements of \erf{two-sys}, are contstructed 
in the second subsection, and the final subsection provides the proof of the Main results.  

\section{A system of algebraic equations and the main results} 
Our main focus is now the system
\[\tag{\ref{two-sys}}\label{eq:2.1}
\bcases
\gl u+f(u-v)=0&\\
\gl v+g(v-u)=0,
\ecases
\]
where $f,g\in C(\R,\R)$ are nondecreasing functions, to be constructed, 
and $\gl>0$ is a constant, to be sent to zero. 

We remark that, due to the monotonicity assumption on $f,g$, the mapping $(u,v)\mapsto (f(u-v),g(v-u)),\,
\R^2\to\R^2$ is monotone. 
Recall that, by definition, a mapping $(u,v)\mapsto (B_1(u,v),B_2(u,v)),\,
\R^2\to\R^2$ is monotone if, whenever $(u_1,v_1), (u_2,v_2)\in\R^2$ satisfy  
$u_1-u_2\geq v_1-v_2$ (resp., $v_1-v_2\geq u_1-u_2$), we have
$B_1(u_1,v_1)\geq B_1(u_2,v_2)$ (resp., $B_2(u_1,v_1)\geq B_2(u_2,v_2)$)

\subsection{Main results}
Our main results are stated as follows. 

\begin{theorem} \label{thm1} There exist two increasing functions $f,g\in C(\R,\R)$ 
having the properties \emph{(a)--(c):} 
\begin{enumerate}
\item[\emph{(a)}]For any $\gl>0$ there exists a unique solution $(u_\gl,v_\gl)\in\R^2$ to 
\erf{two-sys},
\item[\emph{(b)}] the family of the solutions $(u_\gl,v_\gl)$ to \erf{two-sys}, with $\gl >0$, 
is bounded in $\R^2$,
\item[\emph{(c)}] the family $\{(u_\gl,v_\gl)\}_{\gl>0}$ does not converge as $\gl\to 0+$.
\end{enumerate}
\end{theorem} 

It should be noted that, as mentioned in the introduction, the above theorem 
has been somewhat implicitly established by Ziliotto \cite{Zi}. 
In this note, we are interested in a simple and easy approach to finding functions $f,g$ having the properties (a)--(c) in Theorem \ref{thm1}.

The following is an immediate consequence of the above theorem. 

\begin{corollary} \label{cor1}Let $H_i\in C(\R^n,\R)$, $i=1,2$, satisfy $H_1(0)=H_2(0)=0$. Let $f,g\in C(\R,\R)$ be the functions given by Theorem 
\ref{thm1}, and set $B_1(u_1,u_2)=f(u_1-u_2)$ and $B_2(u_1,u_2)=g(u_1-u_1)$ 
for all $(u_1,u_2)\in\R^2$. 
For any $\gl>0$, let $(u_{\gl,1},u_{\gl,2})$ be the (viscosity) solution of \erf{eq:1}. Then, 
the functions $u_{\gl,i}$ are constants, 
the family of the points $(u_{\gl,1},u_{\gl,2})$ 
in $\R^2$ is bounded, and it does not converge as $\gl\to 0+$.   
\end{corollary}

Notice that the convexity of $H_i$ in the above corollary is irrelevant, 
and, for example, one may take $H_i(p)=|p|^2$ for $i\in\I$, which are convex functions.  

We remark that a claim similar to Corollary \ref{cor1} is valid when one 
replaces $H_i(p)$ by degenerate elliptic operators $F_i(x,p,M)$ as far as $F_i(x,0,0)=0$, where $M$ is the variable corresponding to the Hessian matrices of unknown functions.  
(See \cite{CIL} 
for an overview on the viscosity solution approach to 
fully nonlinear degenerate elliptic equations.)

\subsection{The functions $f,g$} 
If $f,g$ are given and $(u,v)\in\R^2$ is a solution of \erf{eq:2.1}, then $w:=u-v$ satisfies
\beq\label{eq:2.2}
\gl w+f(w)-g(-w)=0.
\eeq
Set 
\beq\label{eq:2.3}
h(r)=f(r)-g(-r) \ \ \FOR r\in\R,
\eeq
which defines a continuous and nondecreasing function on $\R$. 

To build a triple of functions $f,g,h$, we need to find two of them in view of the relation \erf{eq:2.3}. We begin by defining function $h$.

For this, we discuss a simple geometry 
on $xy$-plane as depicted in Fig. 1 below. Fix $0<k_1<k_2$. 
The line $y=-\fr 12k_2 +k_1(x+\fr 12)$ has slope $k_1$ and
crosses the lines $x=-1$ and $y=k_2 x$  at 
$\mathrm{P}:=(-1, -\fr 12(k_1+k_2))$ and $\mathrm{Q}:=(-\fr 12,-\fr 12 k_2)$, respectively, while the line $y=k_2x$ meets the lines $x=-1$ and $x=-\fr 12$  at 
$\mathrm{R}:=(-1, -k_2)$ and $\mathrm{Q}=(-\fr 12,-\fr 12 k_2)$, respectively.

Choose $k^*>0$ so that $\fr 12(k_1+k_2)<k^*<k_2$. The line $y=k^*x$ crosses 
the line $y=-\fr 12k_2 +k_1(x+\fr 12)$ at a point $\mathrm{S}:=(x^*,y^*)$ 
in the open line segment between the points $\mathrm{P}=(-\fr 12, -\fr 12(k_1+k_2))$ and 
$\mathrm{Q}=(-\fr 12,-\fr 12 k_2)$. The line connecting $\mathrm{R}=(-1,-k_2)$ and $\mathrm{S}=(x^*,y^*)$ 
can be represented by $y=-k_2+k^+(x+1)$, with $k^+:=\fr{y^*+k_2}{x^*+1}>k_2$.
\begin{figure}
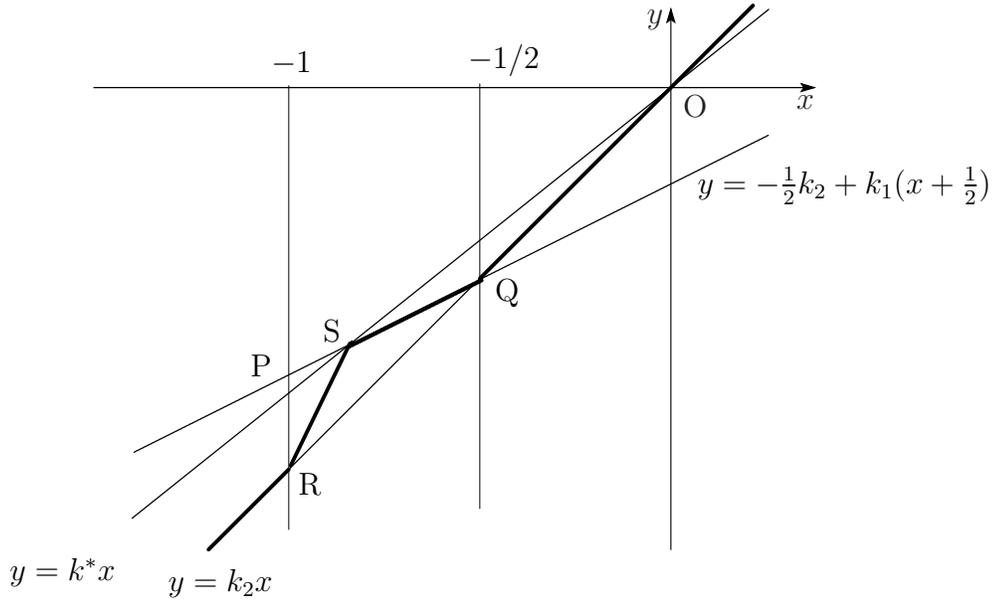
 
\begin{center}
{\input plane-geo1 \hspace{4em}{}}
\caption{Graph of $\psi$.}
\end{center}
\end{figure}

We set 
\[
\psi(x)= \bcases
k_2 x \qquad\qquad \FOR x\in (-\infty, -1]\cup[-1/2,\infty), &\\ 
\min\{-k_2+k^+(x+1), -\fr 12 k_2+k_1(x+\fr 12)\} \ \ \FOR x\in(-1,-\fr 12).&
\ecases
\] 
It is clear that $\psi\in C(\R)$ and increasing on $\R$. 
The building blocks of the graph $y=\psi(x)$ are three lines whose slpoes are 
$k_1<k_2<k^+$.  Hence, if $x_1>x_2$, then $\psi(x_1)-\psi(x_2)\geq k_1(x_1-x_2)$, that is, the function $x\mapsto \psi(x)-k_1 x$ is nondecreasing on $\R$. 

Next, we set for $j\in\N$,
\[
\psi_j(x)=2^{-j}\psi(2^j x)  \ \ \FOR x\in\R.
\]
It is clear that for all $j\in\N$, $\psi_j\in C(\R)$, the function $x\mapsto \psi_j(x)-k_1 x$ is nondecreasing on $\R$, 
and  
\[
\psi_j(x)\bcases
>k_2x \ \ &\FORALL x\in(-2^{-j},-2^{-j-1}), \\
=k_2x \ \ &\OW. 
\ecases
\]
We set
\[
\eta(x)=\max_{j\in\N}\psi_j(x) \ \ \FOR x\in\R.
\]
It is clear that $\eta\in C(\R)$ and $x\mapsto \eta(x)-k_1x$
 is nondecreasing on $\R$. Moreover, we see that 
\[
\eta(x)=k_2 x \ \ \FORALL x\in(-\infty, -\tfrac 12]\cup [0,\infty),
\]
and that if $-2^{-j}< x< -2^{-j-1}$ and $j\in\N$,  
\[
\eta(x)=\psi_j(x)>k_2x. 
\]

Note that the point $\mathrm{S}=(x^*,y^*)$ is on the graph $y=\psi(x)$ and, hence, that 
for any $j\in\N$, the point $(2^{-j}x^*,2^{-j}y^*)$ is on the graph 
$y=\eta(x)$. Similarly, since the point $\mathrm{S}=(x^*,y^*)$ is on the graph $y=k^*x$ and 
for any $j\in\N$, the point $(2^{-j}x^*,2^{-j}y^*)$ is on the graph 
$y=k^*x$. Also, for any $j\in \N$, the point $(-2^{-j},-k_2 2^{-j})$ lies 
on the graphs $y=\eta(x)$ and $y=k_2 x$.   

Fix any $d\geq 1$ and define $h\in C(\R)$ by
\[
h(x)=\eta(x-d). 
\]

For the function $h$ defined above, we consider the problem
\beq\label{eq:2.4}
\gl z+h(z)=0. 
\eeq

\begin{lemma} \label{lem1}For any $\gl\geq 0$, there exists a unique solution $z_\gl\in\R$ of \erf{eq:2.4}. 
\end{lemma}

\bproof Fix $\gl\geq 0$. The function $x\mapsto h(x)+\gl x$ is increasing on $\R$
and satisfies
\[
\lim_{x\to \infty}(h(x)+\gl x)=\infty \ \ \AND
\ \ \lim_{x\to-\infty}(h(x)+\gl x)=-\infty. 
\]
Hence, there is a unique solution of \erf{eq:2.4}. 
\eproof

For any $\gl\geq0$, we denote by $z_\gl$ the unique solution of \erf{eq:2.4}. 
Since $h(d)=0$, it is clear that $z_0=d$.

\begin{lemma} \label{lem2}Let $\gl>0$ and $k>0$. Let $(z,w)\in\R^2$ be the point of the intersection of two lines $y=-\gl x$ and $y=k(x-d)$. Then 
\[
z=\fr{kd}{k+\gl}.
\]
\end{lemma}

\bproof By the assumption, we have
\[
w=-\gl z=k(z-d),
\]
and hence, 
\[
z=\fr{kd}{k+\gl}. \qedhere 
\]
\eproof

\begin{lemma} \label{lem3}There are sequences $\{\mu_j\}$ and $\{\nu_j\}$
of positive numbers converging to zero such that 
\[
z_{\mu_j}=\fr{k_2d}{k_2+\mu_j} \ \ \AND 
\ \ z_{\nu_j}=\fr{k^*d}{k^*+\nu_j}. 
\]
\end{lemma}

\bproof Let $j\in\N$. 
Since $(-2^{-j}, -k_2 2^{-j})$ is on the intersection 
of the graphs $y=k_2x$ and $y=\eta(x)$, it follows that $(-2^{-j}+d, -k_2 2^{-j})$ is on the intersection of the graphs $y=k_2(x-d)$ and $y=h(x)$.  
Set 
\beq \label{eq:2.5}
\mu_j =\fr{k_2 2^{-j}}{d-2^{-j}},
\eeq
and note that $\mu_j>0$ and that
\[
-\mu_j(d-2^{-j})=-k_2 2^{-j},
\]
which says that the point $(d-2^{-j}, -k_22^{-j})$ is on the line 
$y=-\mu_j x$. 
Combining the above with 
\[
-k_2 2^{-j}=h(d-2^{-j})
\]
shows that $d-2^{-j}$ is the unique solution of \erf{eq:2.4}. Also, since 
$(d-2^{-j},-\mu_j(d-2^{-j}))=(d-2^{-j},-k_2 2^{-j})$ is on the line 
$y=k_2(x-d)$, we find by Lemma \ref{lem2} that
\[
z_{\mu_j}=\fr{k_2 d}{k_2+\mu_j}. 
\]

Similarly, since $(2^{-j} x^*, 2^{-j}y^*)$ is on the intersection 
of the graphs $y=k^*x$ and $y=\eta(x)$, we deduce that if we set 
\beq\label{eq:2.6}
\nu_j:=-\fr{2^{-j}y^*}{d+ 2^{-j}x^*}=\fr{2^{-j}|y^*|}{d-2^{-j}|x^*|},
\eeq
then 
\[
z_{\nu_j}=\fr{k^* d}{k^*+\nu_j}. 
\]

It is ovbvious by \erf{eq:2.5} and \erf{eq:2.6} that the sequences 
$\{\mu_j\}_{j\in\N}$ and $\{\nu_j\}_{j\in\N}$ are decreasing and converge to 
zero. 
\eproof 

We fix $k_0\in(0,k_1)$ and define $f, g\in C(\R)$ by $f(x)=k_0(x-d)$ and
\[
g(x)=f(-x)-h(-x). 
\]
It is easily checked that $g(x)-(k_1-k_0)x$ is nondecreasing on $\R$, 
which implies that $g$ is increasing on $\R$, and that 
$h(x)=f(x)-g(-x)$ for all $x\in\R$. We note that
\beq \label{zero of f,g,h}
f(d)=h(d)=g(-d)=0.
\eeq

\subsection{Proof of the main results} We fix $f,g,h$ as above, and consider the system \erf{eq:2.1}. 

\begin{lemma} \label{lem4} Let $\gl>0$. There exists a unique solution of 
\erf{two-sys}. 
\end{lemma}

The validity of the above lemma is well-known, but for the reader's convenience, we provide a proof of the lemma above. 

\bproof By the choice of $f,g$, the functions $f,g$ are nondecreasing on $\R$. We show first the comparison claim: if $(u_1,v_1), (u_2v_2)\in\R^2$ satisfy 
\begin{align}
\label{lem4.1} &\gl u_1+f(u_1-v_1)\leq 0, \quad \gl v_1+g(v_1-u_1)\leq 0,
\\&\gl u_2+f(u_2-v_2)\geq 0, \quad \gl v_2+g(v_2-u_2)\geq 0,
\label{lem4.2}\end{align}
then $u_1\leq u_2$ and $v_1\leq v_2$. Indeed, contrary to this, we suppose 
that $\max\{u_1-u_2,v_1-v_2\}>0$. For instance, if $\max\{u_1-u_2,v_1-v_2\}=u_1-u_2$, then we have $u_1-v_1\geq u_2-v_2$ and $u_1>u_2$, and moreover
\[
0\geq \gl u_1+f(u_1-v_1)\geq\gl u_1+  f(u_2-v_2)>\gl u_2+f(u_2-v_2),
\]
yielding a contradiction. The other case when $\max\{u_1-u_2,v_1-v_2\}=v_1-v_2$, we find a contradiction, $0> \gl v_2+g(v_2-u_2)$, proving the comparison. 

From the comparison claim, the uniqueness of the solutions of \erf{two-sys} 
follows readily. 

Next, we may choose a constant $C>0$ so large that  $(u_1,v_1)=(-C,-C)$ 
and $(u_2,v_2)=(C,C)$ satisfy \erf{lem4.1} and \erf{lem4.2}, respectively. 
We write $S$ for the set of all $(u_1,u_2)\in\R^2$ such that \erf{lem4.1} hold. 
Note that $(-C,-C)\in S$ and that for any $(u,v)\in S$, $u\leq C$ and $v\leq C$.
We set 
\[\bald
u^*&=\sup\{u\mid (u,v)\in S \ \text{ for some }v\},  
\\
v^*&=\sup\{v\mid (u,v)\in S \ \text{ for some }u\}.
\eald \]
It follows that $-C\leq u^*,v^*\leq C$. We can choose sequences 
\[
\{(u_n^1,v_n^1)\}_{n\in\N},\,\{(u_n^2,v_n^2)\}_{n\in\N}\subset S
\] 
such that $\{u_n^1\}, \{v_n^2\}$ are nondecreasing, 
\[
\lim_{n\to\infty}u_n^1=u^* \ \ \AND \ \ 
\lim_{n\to \infty}v_n^2=v^*.
\]
Observe that for all $n\in\N$,  $u_n^2\leq u^*$, $v_n^1\leq v^*$, and 
\[
0\geq\gl u_n^1+f(u_n^1-v_n^1)\geq \gl u_n^1+f(u_n^1-v^*),
\]
which yields, in the limit as $n\to\infty$, 
\[
0\geq \gl u^*+f(u^*-v^*).
\]
Similarly, we obtain $0\geq \gl v^*+g(v^*-u^*)$. Hence, we find that 
$(u^*,v^*)\in S$. 

We claim that $(u^*,v^*)$ is a solution of \erf{two-sys}. Otherwise, we have 
\[
0>\gl u^*+f(u^*-v^*) \ \ \OR \ \ 0>\gl v^*+g(v^*-u^*). 
\]
For instance, if the former of the above inequalities holds, we can choose 
$\ep>0$, by the continuity of $f$, so that 
\[
0>\gl (u^*+\ep)+f(u^*+\ep-v^*). 
\]
Since $(u^*,v^*)\in S$, we have 
\[
0\geq \gl v^*+g(v^*-u^*)\geq \gl v^*+g(v^*-u^*-\ep).
\]
Accordingly, we find that $(u^*+\ep,v^*)\in S$, which contradicts the definition of 
$u^*$. Similarly, if $0>\gl v^*+g(v^*-u^*)$, then we can choose $\gd>0$ so that 
$(u^*,v^*+\gd)\in S$, which is a contradiction. Thus, we conclude that 
$(u^*,v^*)$ is a solution of \erf{two-sys}.   
\eproof

\begin{theorem} \label{thm2}For any $\gl>0$, let $(u_\gl,v_\gl)$ denote the unique solution of \erf{two-sys}. Let $\{\mu_j\}, \{\nu_j\}$ be the sequences 
of positive numbers from Lemma \ref{lem3}. 
Then 
\[
\lim_{j\to\infty} u_{\mu_j}=\fr {k_0d}{k_2}
\ \ \AND \ \ \lim_{j\to\infty} u_{\nu_j}=\fr{k_0d}{k^*}.
\]
In particular,
\[
\liminf_{\gl\to 0}u_\gl\leq \fr{k_0d}{k_2}<\fr{k_0 d}{k^*}\leq \limsup_{\gl \to 0}u_\gl.
\]
\end{theorem}

With our choice of $f,g$, the family of solutions $(u_\gl,v_\gl)$ 
of \erf{eq:2.1}, with $\gl>0$, does not converge as $\gl\to 0$.  

\bproof If we set $z_\gl=u_\gl-v_\gl$, then $z_\gl$ satisfies \erf{eq:2.4}. By Lemma \ref{lem3}, 
we find that 
\[
z_{\mu_j}=\fr{k_2d}{k_2+\mu_j} \ \ \AND \ \ 
z_{\nu_j}=\fr{k^*d}{k^*+\nu_j}. 
\]
Since $u_\gl$ satisfies 
\[
0=\gl u_\gl+f(z_\gl)=\gl u_\gl+k_0(z_\gl-d),
\]
we find that
\[
u_{\mu_j}=-\fr{k_0(z_{\mu_j}-d)}{\mu_j}
=-\fr {k_0d}{\mu_j}\left(\fr{k_2}{k_2+\mu_j}-1\right) 
=-\fr{k_0d}{\mu_j}\fr{-\mu_j}{k_2+\mu_j}
=\fr{k_0d}{k_2+\mu_j}, 
\]
which shows that 
\[
\lim_{j\to\infty}u_{\mu_j}=\fr{k_0d}{k_2}. 
\]
A parallel computation shows that 
\[
\lim_{j\to\infty}u_{\nu_j}=\fr{k_0d}{k^*}.
\]
Recalling that $0<k^*<k_2$, we conclude that  
\[
\liminf_{\gl\to 0}u_\gl\leq \fr{k_0d}{k_2}<\fr{k_0 d}{k^*}\leq \limsup_{\gl \to 0}u_\gl. \qedhere
\]
\eproof 

We remark that, since 
\[
\lim_{\gl \to 0}z_\gl=d \ \ \AND \ \ v_\gl=u_\gl-z_\gl, 
\]
\[
\lim_{j\to\infty}v_{\mu_j}=\fr{k_0d}{k_2}-d \ \ \AND \ \ \lim_{j\to\infty}v_{\nu_j}=\fr{k_0d}{k^*}-d.
\]

We give the proof of Theorem \ref{thm1}. 

\bproof[Proof of Theorem \ref{thm1}] Assertions (a) and (c) are consequences of  Lemma \ref{lem4} and Theorem \ref{thm2}, respectively. 

Recall \erf{zero of f,g,h}. That is,  we have $f(d)=h(d)=g(-d)=0$.  
Setting $(u_2,v_2)=(d,0)$, we compute that for any $\gl>0$,
\[
\gl u_2+f(u_2-v_2)>f(d)=0 
\ \ \AND \ \ \gl v_2 +g(v_2-u_2)=g(-d)=0. 
\]
By the comparison claim, proved in the proof of Lemma \ref{lem4}, we find that 
$u_\gl \leq d$ and $v_\gl\leq 0$ for any $\gl>0$. 
Simlarly, setting $(u_1,v_1)=(0,-d)$, we find that for any $\gl>0$, 
\[
\gl u_1+f(u_1-v_1)=f(d)=0 \ \ \AND \ \ \gl v_1+g(v_1-u_1)\leq g(v_1-u_1)=g(-d)=0,
\]
which shows by the comparison claim that 
$u_\gl \geq 0$ and $v_\gl\geq -d$ for any $\gl>0$. 
Thus, the sequence $\{(u_\gl,v_\gl)\}_{\gl>0}$ is bounded in $\R^2$, which proves assertion (b). 
\eproof

\bproof[Proof of Corollary \ref{cor1}] For any $\gl>0$, we set 
$(u_\gl,v_\gl)\in\R^2$ be the unique solution of \erf{two-sys}. Since 
$H_1(0)=H_2(0)=0$, it is clear that the constant function $(u_{\gl,1}(x),u_{\gl,2}(x)):=(u_\gl,v_\gl)$ is a classical solution of \erf{eq:1}. By a classical uniqueness result (see, for instance, \cite[Theorem 4.7]{IK}), 
$(u_{\gl,1},u_{\gl,2})$ is a unique viscosity solution of \erf{eq:1}. The rest of the claims in Corollary \ref{cor1} is an immediate consequence of Theorem \ref{thm1}.
\eproof

Some remarks are in order. (i) Following \cite{Zi}, we may use Theorem \ref{thm2} as the primary cornerstone for building a scalar Hamilton-Jacobi equation, for which the vanishing discount problem fails to have the full convergence as the discount factor goes to zero.  

(ii) In the construction of the functions $f,g\in C(\R,\R)$ in Theorem \ref{thm2}, the author has chosen $d$ to satisfy $d\geq 1$, but one may choose any 
$d>0$. In the process, the crucial step is to find the function $h(x)=f(x)-g(-x)$, with the properties: (a) the function $x\mapsto h(x)-\ep x$ is nondecreasing on $\R$ for some $\ep>0$ and (b) the curve $y=h(x)$, with $x<d$, meets the lines $y=p(x-d)$ and $y=q(x-d)$, respectively, at $P_n$ and $Q_n$ for all $n\in\N$, where $p,q, d$ are positive constants such that $\ep<p<q$, and the sequences $\{P_n\}_{n\in\N},\, \{Q_n\}_{n\in\N}$ converge to the point $(d,0)$.  Thus, a possible choice of $h$ among many other ways is the following. Define first $\eta\mid \R\to \R$ by $\eta(x)=x(\sin(\log|x|)+2)$ if $x\not=0$, and $\eta(0)=0$ (see Fig. 2). Fix $d>0$ and set $h(x)=\eta(x-d)$ for $x\in\R$.
Note that if $x\not=0$,
\[
\eta'(x)=\sin(\log|x|)+\cos(\log|x|)+2\in[2-\sqrt 2,2+\sqrt 2],
\]
and that if we set $x_n=-\exp(-2\pi n)$ and $\xi_n=-\exp\left(-2\pi n+\fr \pi 2\right)$, $n\in\N$, then 
\[
\eta(x_n)= 2x_n \ \ \AND \ \ \eta(\xi_n)=3\xi_n. 
\]
The points $P_n:=(x_n+d, 2x_n)$ are on the intersection of two curves $y=h(x)$ 
and $y=2(x-d)$, while the points $Q_n:=(d+\xi_n, 3\xi_n)$ are on the intersection of $y=h(x)$ and $y=3(x-d)$. Moreover, $\lim P_n=\lim Q_n=(d,0)$.

\begin{figure} 
\begin{center} 
{
{\unitlength 0.1in%
\begin{picture}(50.0000,34.1000)(1.9000,-36.2000)%
\put(25.9000,-20.1000){\makebox(0,0)[rt]{O}}%
\put(25.6000,-2.1000){\makebox(0,0)[rt]{$y$}}%
\put(51.9000,-20.4000){\makebox(0,0)[rt]{$x$}}%
%
\special{pn 8}%
\special{pa 2600 3620}%
\special{pa 2600 210}%
\special{fp}%
\special{sh 1}%
\special{pa 2600 210}%
\special{pa 2580 277}%
\special{pa 2600 263}%
\special{pa 2620 277}%
\special{pa 2600 210}%
\special{fp}%
%
\special{pn 8}%
\special{pa 190 2000}%
\special{pa 5190 2000}%
\special{fp}%
\special{sh 1}%
\special{pa 5190 2000}%
\special{pa 5123 1980}%
\special{pa 5137 2000}%
\special{pa 5123 2020}%
\special{pa 5190 2000}%
\special{fp}%
\special{pn 13}%
\special{pa 1051 3620}%
\special{pa 1055 3614}%
\special{pa 1060 3606}%
\special{pa 1070 3592}%
\special{pa 1075 3584}%
\special{pa 1090 3563}%
\special{pa 1095 3555}%
\special{pa 1120 3520}%
\special{pa 1125 3514}%
\special{pa 1145 3486}%
\special{pa 1150 3480}%
\special{pa 1155 3473}%
\special{pa 1160 3467}%
\special{pa 1165 3460}%
\special{pa 1170 3454}%
\special{pa 1175 3447}%
\special{pa 1180 3441}%
\special{pa 1185 3434}%
\special{pa 1200 3416}%
\special{pa 1205 3409}%
\special{pa 1235 3373}%
\special{pa 1240 3368}%
\special{pa 1255 3350}%
\special{pa 1260 3345}%
\special{pa 1270 3333}%
\special{pa 1275 3328}%
\special{pa 1280 3322}%
\special{pa 1290 3312}%
\special{pa 1295 3306}%
\special{pa 1305 3296}%
\special{pa 1310 3290}%
\special{pa 1350 3250}%
\special{pa 1355 3246}%
\special{pa 1370 3231}%
\special{pa 1375 3227}%
\special{pa 1380 3222}%
\special{pa 1385 3218}%
\special{pa 1390 3213}%
\special{pa 1395 3209}%
\special{pa 1400 3204}%
\special{pa 1410 3196}%
\special{pa 1415 3191}%
\special{pa 1470 3147}%
\special{pa 1475 3144}%
\special{pa 1485 3136}%
\special{pa 1490 3133}%
\special{pa 1495 3129}%
\special{pa 1500 3126}%
\special{pa 1505 3122}%
\special{pa 1515 3116}%
\special{pa 1520 3112}%
\special{pa 1530 3106}%
\special{pa 1535 3102}%
\special{pa 1575 3078}%
\special{pa 1580 3076}%
\special{pa 1595 3067}%
\special{pa 1600 3065}%
\special{pa 1605 3062}%
\special{pa 1610 3060}%
\special{pa 1615 3057}%
\special{pa 1620 3055}%
\special{pa 1625 3052}%
\special{pa 1630 3050}%
\special{pa 1635 3047}%
\special{pa 1645 3043}%
\special{pa 1650 3040}%
\special{pa 1715 3014}%
\special{pa 1720 3013}%
\special{pa 1730 3009}%
\special{pa 1735 3008}%
\special{pa 1745 3004}%
\special{pa 1750 3003}%
\special{pa 1755 3001}%
\special{pa 1760 3000}%
\special{pa 1765 2998}%
\special{pa 1770 2997}%
\special{pa 1775 2995}%
\special{pa 1785 2993}%
\special{pa 1790 2991}%
\special{pa 1800 2989}%
\special{pa 1805 2987}%
\special{pa 1825 2983}%
\special{pa 1830 2981}%
\special{pa 1925 2962}%
\special{pa 1930 2962}%
\special{pa 2045 2939}%
\special{pa 2050 2937}%
\special{pa 2060 2935}%
\special{pa 2065 2933}%
\special{pa 2075 2931}%
\special{pa 2080 2929}%
\special{pa 2085 2928}%
\special{pa 2095 2924}%
\special{pa 2100 2923}%
\special{pa 2145 2905}%
\special{pa 2150 2902}%
\special{pa 2155 2900}%
\special{pa 2165 2894}%
\special{pa 2170 2892}%
\special{pa 2190 2880}%
\special{pa 2195 2876}%
\special{pa 2205 2870}%
\special{pa 2235 2846}%
\special{pa 2240 2841}%
\special{pa 2245 2837}%
\special{pa 2265 2817}%
\special{pa 2270 2811}%
\special{pa 2275 2806}%
\special{pa 2290 2788}%
\special{pa 2295 2781}%
\special{pa 2300 2775}%
\special{pa 2310 2761}%
\special{pa 2315 2753}%
\special{pa 2320 2746}%
\special{pa 2335 2722}%
\special{pa 2340 2713}%
\special{pa 2345 2705}%
\special{pa 2350 2696}%
\special{pa 2355 2686}%
\special{pa 2360 2677}%
\special{pa 2370 2657}%
\special{pa 2390 2613}%
\special{pa 2405 2577}%
\special{pa 2415 2551}%
\special{pa 2420 2537}%
\special{pa 2425 2524}%
\special{pa 2430 2509}%
\special{pa 2435 2495}%
\special{pa 2445 2465}%
\special{pa 2460 2417}%
\special{pa 2475 2366}%
\special{pa 2495 2294}%
\special{pa 2530 2161}%
\special{pa 2535 2143}%
\special{pa 2540 2124}%
\special{pa 2550 2090}%
\special{pa 2555 2074}%
\special{pa 2560 2059}%
\special{pa 2565 2046}%
\special{pa 2570 2035}%
\special{pa 2575 2026}%
\special{pa 2580 2020}%
\special{pa 2585 2017}%
\special{pa 2595 2013}%
\special{fp}%
\special{pa 2605 1987}%
\special{pa 2615 1983}%
\special{pa 2620 1980}%
\special{pa 2625 1974}%
\special{pa 2630 1965}%
\special{pa 2635 1954}%
\special{pa 2640 1941}%
\special{pa 2645 1926}%
\special{pa 2650 1910}%
\special{pa 2660 1876}%
\special{pa 2665 1857}%
\special{pa 2670 1839}%
\special{pa 2705 1706}%
\special{pa 2725 1634}%
\special{pa 2740 1583}%
\special{pa 2755 1535}%
\special{pa 2765 1505}%
\special{pa 2770 1491}%
\special{pa 2775 1476}%
\special{pa 2780 1463}%
\special{pa 2785 1449}%
\special{pa 2795 1423}%
\special{pa 2810 1387}%
\special{pa 2830 1343}%
\special{pa 2840 1323}%
\special{pa 2845 1314}%
\special{pa 2850 1304}%
\special{pa 2855 1295}%
\special{pa 2860 1287}%
\special{pa 2865 1278}%
\special{pa 2880 1254}%
\special{pa 2885 1247}%
\special{pa 2890 1239}%
\special{pa 2900 1225}%
\special{pa 2905 1219}%
\special{pa 2910 1212}%
\special{pa 2925 1194}%
\special{pa 2930 1189}%
\special{pa 2935 1183}%
\special{pa 2955 1163}%
\special{pa 2960 1159}%
\special{pa 2965 1154}%
\special{pa 2995 1130}%
\special{pa 3005 1124}%
\special{pa 3010 1120}%
\special{pa 3030 1108}%
\special{pa 3035 1106}%
\special{pa 3045 1100}%
\special{pa 3050 1098}%
\special{pa 3055 1095}%
\special{pa 3100 1077}%
\special{pa 3105 1076}%
\special{pa 3115 1072}%
\special{pa 3120 1071}%
\special{pa 3125 1069}%
\special{pa 3135 1067}%
\special{pa 3140 1065}%
\special{pa 3150 1063}%
\special{pa 3155 1061}%
\special{pa 3270 1038}%
\special{pa 3275 1038}%
\special{pa 3370 1019}%
\special{pa 3375 1017}%
\special{pa 3395 1013}%
\special{pa 3400 1011}%
\special{pa 3410 1009}%
\special{pa 3415 1007}%
\special{pa 3425 1005}%
\special{pa 3430 1003}%
\special{pa 3435 1002}%
\special{pa 3440 1000}%
\special{pa 3445 999}%
\special{pa 3450 997}%
\special{pa 3455 996}%
\special{pa 3465 992}%
\special{pa 3470 991}%
\special{pa 3480 987}%
\special{pa 3485 986}%
\special{pa 3550 960}%
\special{pa 3555 957}%
\special{pa 3565 953}%
\special{pa 3570 950}%
\special{pa 3575 948}%
\special{pa 3580 945}%
\special{pa 3585 943}%
\special{pa 3590 940}%
\special{pa 3595 938}%
\special{pa 3600 935}%
\special{pa 3605 933}%
\special{pa 3620 924}%
\special{pa 3625 922}%
\special{pa 3665 898}%
\special{pa 3670 894}%
\special{pa 3680 888}%
\special{pa 3685 884}%
\special{pa 3695 878}%
\special{pa 3700 874}%
\special{pa 3705 871}%
\special{pa 3710 867}%
\special{pa 3715 864}%
\special{pa 3725 856}%
\special{pa 3730 853}%
\special{pa 3785 809}%
\special{pa 3790 804}%
\special{pa 3800 796}%
\special{pa 3805 791}%
\special{pa 3810 787}%
\special{pa 3815 782}%
\special{pa 3820 778}%
\special{pa 3825 773}%
\special{pa 3830 769}%
\special{pa 3845 754}%
\special{pa 3850 750}%
\special{pa 3890 710}%
\special{pa 3895 704}%
\special{pa 3905 694}%
\special{pa 3910 688}%
\special{pa 3920 678}%
\special{pa 3925 672}%
\special{pa 3930 667}%
\special{pa 3940 655}%
\special{pa 3945 650}%
\special{pa 3960 632}%
\special{pa 3965 627}%
\special{pa 3995 591}%
\special{pa 4000 584}%
\special{pa 4015 566}%
\special{pa 4020 559}%
\special{pa 4025 553}%
\special{pa 4030 546}%
\special{pa 4035 540}%
\special{pa 4040 533}%
\special{pa 4045 527}%
\special{pa 4050 520}%
\special{pa 4055 514}%
\special{pa 4075 486}%
\special{pa 4080 480}%
\special{pa 4105 445}%
\special{pa 4110 437}%
\special{pa 4125 416}%
\special{pa 4130 408}%
\special{pa 4140 394}%
\special{pa 4145 386}%
\special{pa 4150 379}%
\special{pa 4155 371}%
\special{pa 4160 364}%
\special{pa 4170 348}%
\special{pa 4175 341}%
\special{pa 4195 309}%
\special{pa 4200 302}%
\special{pa 4215 278}%
\special{pa 4220 269}%
\special{pa 4240 237}%
\special{pa 4245 228}%
\special{pa 4255 212}%
\special{pa 4256 210}%
\special{fp}%
\special{pn 8}%
\special{pa 2060 3620}%
\special{pa 3195 215}%
\special{pa 3197 210}%
\special{fp}%
\special{pn 8}%
\special{pa 980 3620}%
\special{pa 4390 210}%
\special{fp}%
\end{picture}}%
 {}}
\caption{Graph of $\eta$ (slightly deformed).}
\end{center}
\end{figure}
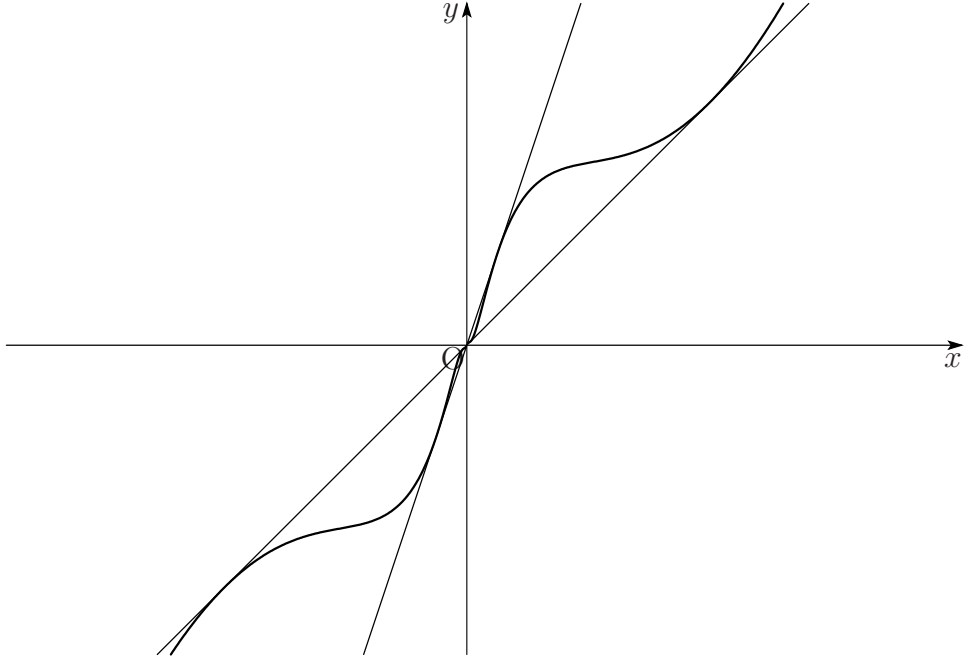 
\section*{Acknowledgments}
The author was supported in part by the JSPS Grants KAKENHI  No. 16H03948, 
No. 20K03688, No. 20H01817, and No. 21H00717.

\bye